\documentclass[10pt,twoside]{article}
\usepackage{amsmath,amssymb,amsthm,amscd}
\usepackage{Latex-document}
\numberwithin{equation}{section}

\def\M{{\mathfrak M}}

\def\cG{{\mathcal G}}

\def\cY{{\mathcal Y}}

\def\bG{{\mathbf G}}

\def\bg{{\mathbf g}}

\def\SO{{\mathbf S}{\mathbf O}}

\def\RR{{\mathbb R}}
\def\CC{{\mathbb C}}

\newtheorem{prop}{Proposition}[section]
\newtheorem{theo}[prop]{Theorem}

\newtheorem{rema}[prop]{Remark}

\def\begeq{\begin{equation}}
\def\endeq{\end{equation}}

\def\and{\quad{\rm and}\quad}

\def\Om{\Omega}

\def\lab{\ }
\def\lab{\label}
\renewcommand{\thesection}{\arabic{section}}
\def\Section#1{\section{\hskip -1em . \hskip 0.6em #1}}
\def\volumeno{I}
\title{\bf Geometry and Nonlinear Analysis\thanks{Supported partially by
NSF grants and a Simons fund.}\vskip 6mm}
\author{Gang Tian\thanks{Department of Mathematics,
Massachusetts Institute of Technology, USA and Beijing University,
China. E-mail: tian@math.mit.edu}\vspace*{-0.5cm}}
\date{\vspace{-8mm}}

\begin{document}

\maketitle \thispagestyle{first} \setcounter{page}{475}

\vskip 12mm

Nonlinear analysis has played a prominent role in the recent
developments in geometry and topology. The study of the Yang-Mills
equation and its cousins gave rise to the Donaldson invariants and
more recently, the Seiberg-Witten invariants. Those invariants
have enabled us to prove a number of striking results for low
dimensional manifolds, particularly, 4-manifolds. The theory of
Gromov-Witten invariants was established by using solutions of the
Cauchy-Riemann equation (cf. [RT], [LT], [FO], [Si], [Ru]). These
solutions are often refered as pseudo-holomorphic maps which are
special minimal surfaces studied long in geometry. It is certainly
not the end of applications of nonlinear partial differential
equations to geometry. In this talk, we will discuss some recent
progress on nonlinear partial differential equations in geometry.
We will be selective, partly because of my own interest and partly
because of recent applications of nonlinear equations. There are
also talks in this ICM to cover some other topics of geometric
analysis by R. Bartnik, B. Andrew, P. Li and X.X. Chen, etc.

Standard partial differential equations in geometry are the
Einstein equation, Yang-Mills equation, minimal surface equation
as well as its close cousin, Harmonic map equation. There are also
parabolic versions of these equations, leading to R. Hamilton's
Ricci flow, the Yang-Mills flow and the mean curvature flow. The
solutions, which played a fundamental role in geometry and
topology, of these equations are their self-dual type ones. I will
focus on self-dual type solutions in this talk. All these
equations are in general hyperbolic equations if we allow Lorentz
metrics on the underlying manifolds, but in differential geometry,
so far, we only concern static solutions, that is, we assume that
the metrics involved are Riemannian. I do believe that the study
of this static case will be very important in our future
understanding general Einstein equation.

\Section{Einstein equation}

\vskip-5mm \hspace{5mm }

We will begin with the Einstein equation. We will always denote by
$M$ a differentiable manifold. A metric $g$ on $M$ is given by a
non-degenerate matrix-valued functions $(g_{ij})$ in local
coordinates $x_1, \cdots, x_n$, where $n$ is the dimension of $M$.
Recall that $g$ is Riemannian if the matrices $(g_{ij})$ are
positive definite.

Associated to each metric, there is a canonical connection, the
Levi-Civita connection, $\nabla$ characterized by the torsion
freeness and $\nabla g=0$, which means that $g$ is parallel. In
local coordinates,
\begin{equation}
\lab{eq:Levi-Civita} \nabla _{\partial\over\partial x_i}
{\partial\over
\partial x_j} = \Gamma^k_{ij} {\partial \over \partial x_k},~~~
\Gamma^k_{ij} = \frac{1}{2} g^{kl} \left ( {\partial g_{il}\over
\partial x_j}+ {\partial g_{jl}\over
\partial x_i}-{\partial g_{ij}\over
\partial x_l}\right ),
\end{equation}
where $(g^{kl})$ denotes the inverse of $(g_{ij})$. Then the
curvature $(R^i_{jkl})$ is defined by
\begin{equation}
\lab{eq:curv} R^i_{jkl} = {\partial \Gamma^i_{jk}\over \partial
x_l} - {\partial \Gamma^i_{jl}\over \partial x_k} + \Gamma^i_{sk}
\Gamma^s_{jl} - \Gamma^i_{sl} \Gamma^s_{jk}.
\end{equation}
The curvature is completely determined by its sectional
curvatures, that is, Gauss curvatures of surface cross sections.
The Ricci curvature ${\rm Ric}=(R_{ij})$ is given by taking trace
of the curvature:
\begin{equation}
\lab{eq:riccicurv} R_{ij} = \sum_k R^k_{ikj}.
\end{equation}
The Ricci curvature essentially measures the variation of the
volume form.

A metric $g$ is called an Einstein metric if it satisfies the
following Einstein equation
\begin{equation}
\lab{eq:einstein} R_{ij} = \lambda g_{ij},
\end{equation}
where $\lambda$ is a constant, usually called Einstein constant.
(\ref{eq:einstein}) is the
Euler-Lagrangian equation of the functional $\int_M s(g) dv$,
where $s(g)$ denotes the scalar curvature of $g$, on the space of
metrics with fixed volume. (\ref{eq:einstein}) is invariant under
diffeomorphism group action. It is elliptic modulo diffeomorphisms
when $g$ is a Riemannian metric. From now on, I will assume that
$g$ is Riemannian.

The simplest examples of Einstein metrics include the euclidean
metric on $\RR^n$, the standard spherical metric on the unit
sphere $S^n$ and the hyperbolic metric on the unit ball $B^n\in
\RR ^n$. In fact, all these metrics have constant sectional
curvature $0$ or $1$ or $-1$, consequently, their Einstein
constant $\lambda$ are $0$, $n-1$, $-n+1$, respectively.

Every Riemannian 2-manifold $(M,g)$ has a natural conformal
structure. The classical Uniformization Theorem states that the
universal covering of $M$ together with the induced conformal
structure is conformal to either $\CC ^1$ or $S^2$ or $B^2$ with
canonical conformal structures. This implies that every Riemannian
2-manifold $(M,g)$ admits a unique Einstein metric within its
associated conformal class and with Einstein constant $0$ or $1$
or $-1$.

Another way of proving this existence is to use partial
differential equation. Given a Riemannian 2-manifold $(M,g)$,
consider a metric $\tilde g$ conformal to $g$, so it is of the
form $e^\varphi g$. A simple computation shows that $\tilde g$ is
of constant curvature $\lambda$ if and only if $\varphi$ satisfies
the following equation
\begin{equation}
\lab{eq:prescribe} \Delta \varphi + {s(g)\over 2} = \lambda
e^\varphi.
\end{equation}
This equation can be solved (cf. \cite{Aubin},
\cite{kazdanwarner}), so there is an Einstein metric in any given
conformal class on any Riemannian 2-manifold. In early 1990's, R.
Hamilton gave a heat flow proof of the Uniformization Theorem
(\cite{Hamiliton}, \cite{Chow}). This new proof also yields a
biproduct: The space of metrics on $S^2$ with positive curvature
is contractible.

The uniformization for 2-manifolds led to two generalizations in
higher dimensions. The one is the Yamabe problem (cf.
\cite{Aubin}, \cite{schoen}). The other is the Calabi's problem on
K\"ahler-Einstein metrics which I will address more later.

For 3-manifolds, Einstein metrics are also of constant sectional
curvature, so their universal coverings are either $S^3$ or
$\RR^3$ or hyperbolic 3-space $H^3$. A major part of Thurston's
program is to show the existence of metrics with constant
sectional curvature on 3-manifolds which satisfy certain
mild topological conditions. Thurston claimed long time ago
that an atoroidal Haken manifold admits a complete hyperbolic
metric. It will be interesting to have an analytic proof of this
claim by solving the corresponding Einstein equation. In general,
one hopes that any 3-manifold can be canonically split into some
pieces of simple topological type and other pieces which admit
Einstein metrics. There are at least two possible approaches to this: one
is the variational method, trying to minimax certain functional
involving curvatures, while the other is to use the Ricci flow,
hoping that one can understand how the singularity is formed along
the flow. So far none of them work yet.

When the dimension is higher than or equal to $4$, an Einstein metric
may not be of constant sectional curvature. It is still a very
interesting question to find out topological constraints on
Einstein manifolds. If a 4-manifold $M$ admits an Einstein metric,
then the Hitchin-Thorpe inequality says that $|\tau(M)| \le
\frac{2}{3} |\chi(M)|$. This implies that the connected sum of
more than $4$ copies of $\CC P^2$ can not have any Einstein
metric. More recently, C. Lebrun showed that a 4-manifold $M$ with
non-vanishing Seiberg-Witten invariant admits an Einstein metric
only if $3 \tau(M) \le \chi(M)$. We do not know any constraints on
compact Einstein manifolds of dimension higher than $4$. It may be
possible that any manifold of dimension $\ge 5$ has an Einstein
metric.

Examples of Einstein metrics can be constructed by exploring
symmetries, such as, homogeneous Einstein metrics, cohomogeneity
one Einstein metrics. One can also construct new
Einstein metrics from known ones through certain intrigue
constructions when the underlying manifolds are of special
fibration structures  (cf. \cite{wangmc}, \cite{boyergalicki}).

When $n\ge 4$, there is a special class of solutions of the
Einstein equation, that is, Einstein metrics of special holonomy.
If $(M,g)$ is an irreducible Riemannian manifold, a well-known
theorem of M. Berger states that either $(M,g)$ is a locally
symmetric space or its reduced holonomy is one of the following
groups: $SO(n)$, $U(\frac{n}{2})$ ($n\ge 4$), $SU(\frac{n}{2})$ ($n\ge 4$),
$Sp(1)\cdot Sp(\frac{n}{4})$ ($n\ge 8$), $Sp(\frac{n}{4})$ ($n\ge 4$)
and two exceptional groups
$Spin(7)$ and $G_2$. We call $(M,g)$ a Riemannian manifold with
special holonomy if it is irreducible and its (reduced) holonomy is strictly
contained in $SO(n)$. It can be shown that a Riemannian manifold
of special holonomy is automatically Einstein if its holonomy is
other than $U(\frac{n}{2})$. We have K\"ahler-Einstein metrics for
the $U(\frac{n}{2})$ case. In fact, these special Einstein metrics
are of self-dual type. Each manifold $(M,g)$ of special holonomy
has a parallel $n-4$ form defined as follows: Let $W\subset
\Lambda^2TM$ be the subspace generated by
the Lie algebra of the holonomy group of $(M,g)$, then the curvatures lie in $S^2(W)$.
Define a $4$-form $\psi(W)$ by
\begin{equation}
\lab{eq:4-form} \psi(W)=\sum w_i\wedge w_i,
\end{equation}
where $\{w_i\}$ is an orthonormal basis of $W$. Clearly, it is
independent of the choice of $\{w_i\}$ and is parallel. This
4-form induces a symmetric operator $T_{\psi(W)}: W\mapsto W$:
$T_{\psi(W)}(v) = i_v\psi(W)$, where $i_v$ denotes the interior
product with $v$. One can check that $T_{\psi(W)}$ has at most two
distinct eigenvalues. Moreover, there is a distinguished
eigenspace $W_0$ of $T_{\psi(W)}$ of codimension $0$, $1$ and $3$.
Let $\beta$ be the corresponding eigenvalue. Put
\begin{equation}
\lab{eq:omega} \Omega(W)= \frac{1}{\beta} \ast \psi(W).
\end{equation}
Clearly, it is parallel. Denote by $S^2W = S^2W_0 + S^2W_1$ the
decomposition according to eigenvalues, then the curvature $R(g)$
of $g$, which lies in $S^2W$, is decomposed into $R_0\in S^2W_0$
and $R_1\in S^2W_1$. Furthermore, we have
\begin{equation}
\lab{eq:selfdual} R_0\wedge \Omega = \ast R_0
\end{equation}
and $R_1$, which can be void, is completely determined by Ricci
curvature of $g$. Therefore, manifolds of special holonomy are
always self-dual. It is very important in the study of Einstein
metrics of special holonomy. For example, the self-duality implies
an a prior $L^2$-bound on curvature: There is a uniform constant
$C(p_1(M), \Omega, s(g))$, depending only on the first Pontrjagin
class, $\Omega$ and the scalar curvature $s(g)$, such that for any
Einstein metric $g$ of special holonomy, we have
\begin{equation}
\lab{eq:l2curvature} \int_M |R(g)|^2 dv = C(p_1(M), \Omega, s(g)).
\end{equation}
Here $\Omega$ is the corresponding parallel form. In dimension $4$,
we can study self-dual Einstein metrics, that is, Einstein metrics with self-dual
Weyl curvature. These self-dual metrics share similar
properties as those with special holonomy do.

The special geometry we see most is the K\"ahler geometry. A
K\"ahler manifold is a Riemannian manifold $(M,g)$ whose holonomy
lies in $U(\frac{n}{2})$, it is equivalent to saying that $M$ has a
compatible and parallel complex structure $J$, that is, $g(Ju,
Jv)=g(u,v)$, where $u, v \in TM$ are arbitrary, and $\nabla J=0$.
So $M$ is a complex manifold with induced complex structure
by $J$. Usually, we denote $g$ by its K\"ahler form $\omega_g =
g(\cdot, J\cdot)$. In local complex coordinates $z_1,\cdots, z_m$ of $M$
($n=2m$),
\begin{equation}
\lab{eq:kahlerform} \omega_g = \frac{\sqrt{-1}}{2} \sum_{i,j=1}^m
g_{i\bar j} dz_i \wedge d\overline z_j ,
\end{equation}
where $(g_{i\bar j})$ is a positive Hermitian matrix-valued
function. The self-duality simply means that the curvature of a
K\"ahler metric has only components of type (1,1). A K\"ahler
metric $g$ is Einstein if and only if the trace of its curvature
against $\omega_g$ is constant, we call such a metric
K\"ahler-Einstein. A necessary condition for the existence of
K\"ahler-Einstein metric on $M$ is that the first Chern class
$c_1(M)$ is definite. Since the Ricci curvature of a K\"ahler
metric $g$ can be expressed as $-\partial\overline\partial \log
det(g_{i\bar j})$, the Einstein equation is reduced to solving the
following complex Monge-Amper\'e equation
\begin{equation}
\lab{eq:mongeampere} \det(g_{i\bar j} + \frac{\partial ^2
\varphi}{\partial z_i\partial {\overline z}_j}) = e^{h- \lambda
\varphi} \det(g_{i\bar j}),~~~~(g_{i\bar j} + \frac{\partial ^2
\varphi}{\partial z_i\partial {\overline z}_j}) > 0,
\end{equation}
where $\varphi$ is unknown and $h$ is a given function depending
only on $g$. This is a fully nonlinear elliptic equation and
easier to solve.

A program initiated by E. Calabi in early 1950's is to study the
existence and uniqueness of K\"ahler-Einstein
metrics.\footnote{Later in 1980's, E. Calabi extended this to
extremal K\"ahler metrics, one can see X.X. Chen's paper in this
proceeding for recent progresses on extremal metrics.} The
uniqueness of K\"ahler-Einstein metrics was known in 1950's in the
case that the first Chern class is nonpositive and was proved by
Bando-Mabuchi \cite{bandomabuchi} in 1986 in the case that the
first Chern class is positive. The difficult part of Calabi's
program is about the existence. The celebrated solution of Yau
[Ya] for the Calabi conjecture established the existence of a
Ricci-flat metric, now named as Calabi-Yau metric, in each
K\"ahler class on a compact K\"ahler manifold $M$ with $c_1(M)=0$.
If $c_1(M) < 0$, the existence of K\"ahler-Einstein metrics was
proved by Yau [Ya] and Aubin [Au], independently. There are
further analytic obstructions to the existence of
K\"ahler-Einstein metrics on $M$ with $c_1(M)
>0$. Matsushima proved that $M$ has a K\"ahler-Einstein metric only
if the Lie algebra $\eta (M)$ of its holomorphic fields is
reductive. Also if $M$ has a K\"ahler-Einstein metric, then the
Futaki invariant from \cite{futaki} vanishes. The Futaki invariant
is a character of $\eta (M)$. If $M$ is a complex surface with
$c_1(M)>0$, then it admits a K\"ahler-Einstein metric if and only
if the Lie algebra of holomorphic vector fields is reductive
\cite{tian1}. For a general $M$ with $c_1(M)>0$, the existence of
K\"ahler-Einstein metrics is equivalent to certain analytic
stability \cite{tian2}. This analytic stability amounts to
checking an nonlinear inequality of Moser-Trudinger type: Assume
that $\eta(M)=\{0\}$. \footnote{If $\eta(M)\not= \{0\}$, then the
inequality holds only for those functions perpendicular to
functions induced by holomorphic vector fields.} If $\omega$ is a
K\"ahler metric with $[\omega]=c_1(M)$ and $\varphi$ with $\int_M
\varphi \omega^n =0$ and $\omega
+\partial\overline{\partial}\varphi >0$,
\begin{equation}
\lab{eq:mosertrudinger}
\log\left (\int_M e^{-\varphi} \omega^n\right )
\le J_\omega(\varphi) - f(J_\omega(\varphi)),
\end{equation}
where $f$ is some function bounded from below and satisfies
$\lim_{t\to \infty}f(t) = \infty$\footnote{$f$ may depend on $\omega$.} and $J_\omega$
is defined by
\begin{equation}
\lab{eq:Jenergy}
J_\omega(\varphi) = \sum_{i=0}^{n-1} {i+1\over n+1} {\sqrt{-1}\over 2 V} \int_M \partial \varphi
\wedge {\overline{\partial}}\varphi \omega^i \wedge
(\omega +\partial\overline{\partial}\varphi )^{n-i-1},
\end{equation}
where $V = \int_M \omega^n$. The inequality (\ref{eq:mosertrudinger})
has been checked for many manifolds, such as Fermat hypersurfaces.
Furthermore, the analytic stability implies the asymptotic CM-stability of $M$
introduced in \cite{tian2} in terms of Geometric Invariant Theory. If one proved the
partial $C^0$-estimate conjectured in \cite{tian3}, this asymptotic
stability in [Ti2] would imply the existence of K\"ahler-Einstein
metrics. Very recently, by using the Tian-Yau-Zeldich expansion
(cf. \cite{tian4}, \cite{catlin}, \cite{zelditch})
and a result of Zhiqin Lu \cite{zhiqin}, S. Donaldson \cite{donaldson1}
proved the asymptotic Chow
stability \cite{mumford} of algebraic manifolds which admit
K\"ahler-Einstein metrics \cite{donaldson1}. This gives a partial
answer to one conjecture of Yau: If $\eta(M)=0$, then there is a
K\"ahler-Einstein metric on $M$ if and only if $M$ is asymptotically
Chow stable. It would be a very interesting problem in algebraic
geometry to compare the Chow stability with the CM-stability introduced in
[Ti2]. Both stabilities can be defined in terms of the Chow coordinate
of $M$, but their corrsponding polarizations are different
(cf. \cite{paultian}).

K\"ahler-Ricci solitons arose naturally from the study of
the existence of K\"ahler-Einstein metrics and
Hamilton's Ricci flow in K\"ahler geometry and generalize
K\"ahler-Einstein metrics. A K\"ahler metric $g$ is a
K\"ahler-Ricci soliton if there is a holomorphic field $X$ such
that
\begin{equation} \lab{eq:soliton} {\rm Ric}(g) -\lambda
\omega_g = L_X\omega_g.
\end{equation}
As before, this equation can be reduced to a sightly more complicated
complex Monge-Ampere equation (cf.\cite{tianzhu}).
It was proved in \cite{tianzhu} that K\"ahler-Ricci solitons are unique
modulo automorphisms. In subsequent papers, we also gave an analytic
criterion for the existence as one did in [Ti2]. It was conjectured that
given any K\"ahler manifold $M$ with $c_1(M) >0$, either $M$ has a
K\"ahler-Einstein metric or there are diffeomorphisms $\phi_i$ and
K\"ahler metrics $g_i$ such that $\phi_i^*g_i$ converge to a
unique K\"ahler-Ricci soliton on $M'$, which may be different from
$M$. This conjecture was posed by R. Hamilton in studying the
Ricci flow and myself in studying K\"ahler-Einstein metrics.
When the complex dimension of $M$ is $2$, in view of the main result in \cite{tian1},
it suffices to show that the blow-up of $\CC P^2$ at two points admits a K\"ahler-Ricci
soliton. This should be doable.

So far, the most successful method in proving the existence is the
continuity method. The other possible approach is to use the
K\"ahler-Ricci flow, which has only partial success (cf.
X.X.Chen's talk at this ICM).

There remain many problems in studying K\"ahler-Einstein with
prescribed singularities, though a lot has been done (cf.
\cite{chengyau}, \cite{tianyau1}, \cite{tsuji}, etc.). A given
K\"ahler manifold $M$ may not have definite first Chern class, so
it does not admit any K\"ahler-Einstein metrics, but by blowing
down certain subvarieties, the resulting manifold (possibly
singular) may admit a canonical K\"ahler-Einstein metric. For
instance, if $M$ is an algebraic manifold of general type, can any
given K\"ahler metric be deformed along the K\"ahler-Ricci flow to
a unique K\"ahler-Einstein metric? Is the limiting metric
independent of the initial metric? The answer to these questions
seems to be affirmative in complex dimension $2$ or in the case
that minimal models exist. Another unsolved problem is Yau's
conjecture: Every complete Calabi-Yau open manifold can be
compactified such that the divisor at infinity is the zero-locus
of a section of a line bundle proportional to the anti-canonical
bundle. This is a hard problem. In \cite{tianyau2},
\cite{tianyau3} and \cite{bandokobayashi}, complete Calabi-Yau
manifolds were constructed on complements of a smooth divisor
which is a fraction of anti-canonical divisor and satisfies
certain positivity conditions (also see \cite{joyce}). In view of
these and \cite{cheegertian1}, one is led to the following
conjecture: a complete Calabi-Yau manifold $M$ with quadratic
curvature decay and euclidean volume growth is of the form
$M={\overline M}\backslash D$ such that $D$ is ample near $D$ and
the anti-canonical bundle $K_{\overline M}^{-1}$ is $ \alpha [D]$
for some $\alpha > 1$. This can be considered as the refinement of
Yau's conjecture in a special case.

The next special holonomy is contained in $Sp(1)Sp(\frac{n}{4})$.
Riemannian manifolds with such a holonomy are called
quaternion-K\"ahler manifolds. They are automatically Einstein.
The prototype is the quaternionic projective space ${\bold P}^
{\frac{n}{4}}_{\bold H}$. There are many examples of
quaternion-K\"ahler manifolds due to the works of many people,
including S. Salamon, Galicki-Lawson, Lebrun, etc.
Quaternion-K\"ahler manifolds with zero scalar curvature are
hyperK\"ahler, that is, its holonomy lies in $Sp(\frac{n}{4})$.
The existence of hyperK\"ahler manifolds follows from Yau's
solution for the Calabi conjecture. However, we do not know yet if
there are quaternion-K\"ahler manifolds with positive scalar
curvature and which are not locally symmetric, while we do have a
number of symmetric ones, the so-called Wolf spaces. It led Lebrun
and S. Salamon to guess that the Wolf spaces are all complete
quaternion-K\"ahler manifolds with positive scalar curvature. So
far, it has been checked up to dimension $\frac{n}{4}\le 3$. We
would like to point out that there are many non-symmetric
quaternion-K\"ahler orbifolds with positive scalar curvature due
to Galicki-Lawson (\cite{galickilawson}).

Riemannian manifolds with holonomy $G_2$ and ${\rm Spin}(7)$ must
be Ricci-flat and of dimension $7$ and $8$, respectively. It took
a long time to settle the question of whether such metrics exist,
even locally. Local metrics with these holonomy were constructed
by R. Bryant \cite{bryant}. Later, complete examples were
constructed by Bryant and S. M. Salamon\cite{bryantsalamon}.
Examples of compact 7- and 8-manifolds with holonomy $G_2$ and
${\rm Spin}(7)$ were first constructed by D. Joyce in early 1990's
(cf. \cite{joyce}). D. Joyce's construction was inspired by the
Kummer construction: metrics with holonomy ${\rm SU}(2)$ on the
$K3$ surface can be obtained by resolving the $16$ singularities
of the orbifolds $T\sp 4 /Z\sb 2$, where $Z\sb 2$ acts on $T\sp 4$
with $16$ fixed points. In the case of $G\sb 2$, Joyce chooses a
finite group $\Gamma\subset G_2$ of automorphisms of the torus
$T\sp 7$. Then he resolves the singularities of $T\sp 7/\Gamma$ to
get a compact $7$-manifold $M$ with holonomy $G_2$. A similar
construction can be implemented for the ${\rm Spin}(7)$ case by
choosing a finite group $\Gamma$ of automorphisms of the torus
$T\sp 8$ and a flat $\Gamma$-invariant ${\rm Spin}(7)$-structure
on $T\sp 8$. More recently, Kovalev gave a new construction of
Riemannian metrics with special holonomy $G_2$ on compact
7-manifolds. The construction is based on gluing asymptotically
cylindrical Calabi-Yau manifolds built up on the work in
\cite{tianyau2}. Examples of new topological types of compact
7-manifolds with holonomy $G_2$ were obtained.

So far, all Ricci-flat compact manifolds are of special holonomy.
There should exist complete Ricci-flat manifolds with generic
holonomy ${\rm SO}(n)$. The question is how we can find them. Here
is a possible example in $4$-dimension: It was shown that there is
a Calabi-Yau manifold with cylindrical end asymptotic to $T^3$ \cite{tianyau2}.
Complex analytically, this manifold can be obtained by blowing up
the $9$ base points of a generic elliptic pencil on $\CC P^2$ and
removing one smooth fiber. Now take two copies of such Calabi-Yau
manifolds and glue them along the $T^3$'s at infinity. One way of
gluing them is to respect the complex structures, then we will get
a $K3$ surface. Could one use different gluing maps which do not
preserve the complex structures, so that one may obtain new
Ricci-flat manifolds with generic holonomy? We can also ask if any
complete Ricci-flat manifolds can be decomposed in some sense into
a connected sum of Calabi-Yau manifolds.
Similar things can be done in higher dimensions.

Geometry of moduli space of Einstein manifolds is extremely
important. For example, the moduli space of Calabi-Yau manifolds
provides the B-model in the Mirror Symmetry. If $(M,g)$ is an
Einstein manifold with special holonomy, then it was proved that
nearby Einstein manifolds in the moduli space is also of special
holonomy. The first analytic problem about the
moduli is its compactness. The moduli space is very often
noncompact, so we need to compactify it. Then we can consider what
structures a compactified moduli space has.

We have pointed out before that for any Einstein manifold
$(M,g)$ with special holonomy, the $L^2$-norm of its curvature
depends only on the second Chern character and the Einstein
constant $\lambda$ (assuming that the volume of $M$ is normalized,
say $1$). One can first give a weak compactification $\overline
{\mathcal M}$ of the moduli space ${\mathcal M}$ of Einstein
manifolds with special holonomy in the Gromov-Hausdorff topology.
A basic problem is the regularity of a limit in ${\overline
{\mathcal M}}\backslash {\mathcal M}$. There are two cases of the
limit, one is when the limit is still compact, while the other has
infinite diameter as a length space. Here let us consider only the
first case, since we know much more in this case and it is
necessary for studying the second. If $M_\infty$ is a compact
limit, then there is a sequence of Einstein manifolds $(M_i, g_i)$
with special holonomy and bounded diameter converging to
$M_\infty$ in the Gromov-Hausdorff topology. When the dimension is
$2$, it was proved in \cite{tian1} that $M_\infty$ is a
K\"ahler-Einstein orbifold with isolated singularities. Its real
version was done by M. Anderson in \cite{anderson}. The
compactness theorem played a very important role in the resolution
of the Calabi problem for complex surfaces (cf. \cite{tian1}). In \cite{cheegertian2},
Cheeger and I proved

\begin{theo}\footnote{The K\"ahler case of this theorem was proved
in \cite{cheegercoldingtian}.}
\lab{th:cheegertian} Let $M_\infty$ be
the above limit of a sequence of Einstein manifolds $(M_i,g_i)$
with the same special holonomy and uniformly bounded diameter.
Then there is a rectifiable closed subset $S\subset M_\infty$ such
that $M_\infty\backslash S$ is a smooth manifold which admits an
Einstein metric $g_\infty$ with the same holonomy as $(M_i,g_i)$
do. Furthermore, $M_\infty$ is the metric completion of
$M_\infty\backslash S$ with respect to the distance induced by
$g_\infty$.
\end{theo}

This is based on deep works of Cheeger-Colding \cite{cheegercolding}
on spaces which are limits of
manifolds with Ricci curvature bounded from below,
my joint work with Cheeger and Colding on structure of the singular sets of
limits of Einstein manifolds with $L^2$ curvature bounds
\cite{cheegercoldingtian} and
Cheeger's work on rectifiability of singular sets of the limits
\cite{cheeger}).

\begin{rema}
\lab{rema: cheegertian} {\rm In fact, the convergence can also be
strengthened: There is an exhausion of $M_\infty\backslash S$ by
compact sets $K_i\subset K_{i+1}\cdots$ and diffeomorphisms
$\phi_i: K_i\mapsto M_i$ such that $\phi_i(K_i)$ converge to $S$
in the Gromov-Hausdorff topology and $\phi^*g_i$ converge to
$g_\infty$ in the $C^\infty$-topology.}
\end{rema}

The most fundamental problem left is the regularity of $S$ or
structure of $M_\infty$ along $S$. The conjecture is that $S$ can
be stratified into $\coprod_{i\le n-4} S_i$ such that each stratum
$S_i$ is a smooth manifold of dimension $\le i$. If $(M_i, g_i)$ are
K\"ahler-Einstein, then $S_{2j+1}=S_{2j}$. We know (cf.
\cite{cheegercoldingtian},
\cite{cheegertian2}) that tangent cones at almost all points of
$S_{n-4}$ are of the form $\RR^{n-4}\times C(S^3/\Gamma)$, where
$C(S^3/\Gamma)$ is the cone over $S^3/\Gamma$ and $\Gamma \subset
SO(4)$ is a finite group. It makes us conjecture that $M_\infty$
should be homeomorphic to an open set of $\RR^{n-4}\times C(S^3/\Gamma)$
locally along $S_{n-4}\backslash\coprod_{i\le n-5} S_i$. It is plausible that
$M_\infty$ is actually smooth along $S_{n-4}\backslash\coprod_{i\le n-5} S_i$
in a suitable sense. We also believe that the $(n-4)$-form $\Omega_\infty$
associated to the special holonomy of $g_\infty$ extends to $S$ in a suitable
sense and its restriction to $S$ is the same as the volume, i.e.,
$S$ is calibrated by $\Omega_\infty$ in a suitable sense.

When $(M_i,g_i)$ are K\"ahler-Einstein manifolds with positive scalar curvature,
it was conjectured by the author long time ago that a multiple of the anticanonical bundle of
$M_\infty\backslash S$ extends to be a line bundle across the singular set $S$.
This is of course true if one can show that $M_\infty$ has only quotient
singularities. The affirmation of this conjecture will enable
us to prove the converse of a result in \cite{tian2}, that is,
the algebraic stability of a K\"ahler manifold with
positive first Chern class assures the existence of K\"ahler-Einstein metrics.

Finally, we shall refer the readers to \cite{cheegertian2} for detailed study of
tangent cones at any singularity of $M_\infty$.

\Section{Yang-Mills equation}

\vskip-5mm \hspace{5mm }

Next we discuss the Yang-Mills equation. The Yang-Mills equation
has played a fundamental role in our study of physics and geometry
and topology in last few decades. In the following, unless
specified, we assume that $(M,g)$ is a Riemannian manifold of
dimension $n$ and $\bG$ is a compact subgroup in $\SO(r)$ and
$\bg$ is its Lie algebra. Let $E$ be a $\bG$-bundle over $M$.

First we recall that a connection of $E$ over $M$ is locally of
the form
\begin{equation}
\lab{eq:1.1}
A = A_i dx_i, \quad A_i \in \bg
\end{equation}
where $x_1, \cdots , x_n$ are euclidean coordinates of $\RR^n$ and
$A_i$
are matrices in $\bg$.
Its curvature can be computed as follows:
\begin{equation}
\lab{eq:1.2} F_A = dA + A\wedge A .
\end{equation}

The Yang-Mills functional is defined on the space of connections and
given by
%%%%%%%%%%%%%%%
\begin{equation}
\label{eq:1.4} \cY(A)=\frac{1}{4\pi^2}\int_M |F_A|^2_g dV_g.
\end{equation}
%%%%%%%%%%%%%%%
The Yang-Mills equation is simply its Euler-Lagrange equation
\begin{equation}
\lab{eq:yangmills} D_A^*F_A=0,
\end{equation}
where $D_A$ denotes the covariant derivative of $A$ and $D_A^*$ is
its adjoint. On the other hand, being the curvature of a
connection, $A$ automatically satisfies the second Bianchi
identity $D_AF_A=0$. We will call $A$ a Yang-Mills connection if
it satisfies (\ref{eq:yangmills}).

The gauge group $\cG$ consists of all sections of $Ad(E)$ over
$M$, locally, they are just maps into $\bG\subset \SO(r)$. It acts
on the space of connections by assigning $A$ to $\sigma(A)=\sigma
A\sigma^{-1} - \sigma d \sigma^{-1}$ for each $\sigma \in \cG$.
Clearly, the Yang-Mills functional is invariant under the action
of $\cG$, so does the Yang-Mills equation. In particular, it
implies that the Yang-Mills equation is not elliptic. A difficult problem is
to construct good gauges which can be controled by curvatures.
So called Coulomb gauges have been constructed by Uhlenbeck \cite{uhlenbeck2}
in $L^{n/2}$-norms and more recently, by Tao-Tian and Meyer-Riviere in
Morrey norms (cf. \cite{taotian}).

The simplest Yang-Mills connections are provided by harmonic one
forms: If $G=U(1)$, then $\bg= i \RR$ and $A$ is simply a one-form
and the Yang-Mills equation is $d^* d A=0$, the gauge
transformation is given by $\sigma= e^{i a}\mapsto A + i d a$. It
follows that modulo gauge transformations, abelian Yang-Mills
connections are in one-to-one correspondence with harmonic one
forms.

Now we assume that $(M,g)$ is of special holonomy. Let $\Omega$ be
the associated closed form of degree $n-4$. We say that a
connection $A$ is $\Omega$-self-dual if
\begin{equation}
\lab{omegaselfdual} \ast(\Omega\wedge F_A) =  F_A,
\end{equation}
where $\ast$ is the Hodge operator \footnote{If $\bG\subset U(r)$,
one can consider more general self-dual equation: $A$ is an $\Omega$-self-dual
connection if ${\rm tr}(F_A)$ is harmonic
and $\ast(\Omega\wedge F_A^0) =  F_A^0$, where $F^0_A=F_A - \frac{1}{r} {\rm tr}(F_A) Id$.
If $\bG={\rm SU}(r)$, this coincides with (\ref{omegaselfdual}).}

One can show that an $\Omega$-self-dual connection is a Yang-Mills
connection. Clearly, the self-duality is invariant under gauge
transformations. So self-dual connections provide a special class
of Yang-Mills solutions.

There are many examples of $\Omega$-self-dual connections. First,
the Levi-Civita connection of the underlying Riemannian metric is
$\Omega$-self-dual. In this sense, the Yang-Mills equation is a
semi-linear version of the Einstein equation. Secondly, if $E$ is
a stable holomorphic vector bundle, then the
Donaldson-Yau-Uhlenbeck theorem (\cite{donaldson2},
\cite{yauuhlenbeck}) states that $E$ has a unique
Hermitian-Yang-Mills connection, an easy computation shows that a
connection is Hermitian-Yang-Mills if and only if it is
$\Omega$-self-dual, where $\Omega=
-\frac{\omega^{n/2-2}}{(n/2-2)!}$ and $\omega$ is the underlying
K\"ahler form. Thirdly, if $(M,g)$ is a Calabi-Yau 4-fold and
$\Omega$ is its associated $(n-4)$-form induced by the
$SU(4)\subset Spin(7)$- structure, then $\Omega$-self-dual
connections are just complex self-dual instantons of
Donaldson-Thomas \cite{donaldsonthomas}. Also, one may construct
$\Omega$-self-dual instantons from $\Omega$-calibrated
submanifolds of $M$.

When $n=4$, self-dual instantons were used to construct the
Donaldson invariants for 4-manifolds. This eventually led to the
Seiberg-Witten invariants, which is much easier to compute. The
construction goes roughly as follows: Let $M$ be a $4$-manifold
and $g$ is a generic metric. Let $E$ be an $SU(2)$-bundle over
$M$. Consider the moduli space ${\mathcal M}_E$ of self-dual
instantons of $E$, that is, solutions of
\begin{equation}
\lab{eq:antiselfdual} F_A = \ast F_A
\end{equation}
modulo gauge transformations. A generalized instanton consists of
an anti-self-dual instanton and a tuple of points of $M$ such that
the second Chern class of the instanton and totality of the points
sum up to represent $c_2(E)$. Let $\overline{\mathcal M}_E$ be the
moduli space of all generalized instantons of $E$. Then Uhlenbeck
compactness theorem states that $\overline{\mathcal M}_E$ is
compact. Also $\overline{\mathcal M}_E$ is a stratified space with
${\mathcal M}_E$ as its main stratum. If $b_2^+(M) \ge 3$ and $g$
is generic, then each stratum has expected dimension which can be
easily computed by the Atiyah-Singer index theorem, so
$\overline{\mathcal M}_E$ can be taken as a fundamental class. The
Donaldson invariants are obtained by integrating pull-backs of
cohomology classes of $M$ on this fundamental class.

Similarly, one can define the Seiberg-Witten invariant by using
the Seiberg-Witten equation. Technically, it is much easier since
the moduli space is already compact. Sometimes, it was said that
the Seiberg-Witten invariant does not need hard analysis, in fact,
it is false. Taubes' deep theorem on equivalence of Seiberg-Witten
and Gromov-Witten invariants requires hard analysis.

What about higher dimensional cases? Can we construct new
deformation invariants by using $\Omega$-self-dual
connections? In order to achieve it, one has to consider the
following issues: 1. Is the corresponding self-dual equation
elliptic? Indeed, the self-dual equation on a $Spin(7)$-manifold
is elliptic. 2. Can we compactify the moduli space? If so, how do
we stratify the compactified moduli space? 3. Does each stratum
have right dimension which can be predicted by the index
theorem? If we solve these issues, we can define new invariants,
then we can study how to compute them.

The first issue is easy to check. We just need to linearize the
self-dual equation and see if it is elliptic. There are examples,
such as, self-dual equations on 4-manifolds and Donaldson-Thomas
complex self-dual on Calabi-Yau 4-manifolds. It will be very
useful to construct deformation invariants by using complex
self-dual instantons. The success of it will provide a powerful
tool of constructing holomorphic cycles of codimension $4$, which
are pretty much evading us.

Next we consider the compactification. Having a good compactification,
we will be able to get property 3 in the above.
Let $(M,g)$ be a compact
Riemannian manifold of dimension $n$ and with special holonomy.
Let $\Omega$ be the associated closed form of degree $n-4$. Let
$E$ be a unitary vector bundle over $M$. Recall that
${\M}_{\Omega, E}$ consists of all gauge equivalence classes of
$\Omega$-asd instantons of $E$ over $M$. In general,
${\M}_{\Omega, E}$ may not be compact. So we will compactify it.

An admissible $\Omega$-self-dual instanton is simply a smooth connection $A$ of $E$
over $M\backslash S(A)$ for a closed subset
$S(A)$ of Hausdorff dimension $n-4$ such that $\int_M |F_A|^2 < \infty$.
A generalized $\Omega$-self-dual instanton is made of an
admissible $\Omega$-self-dual instanton $A$ of $E$ and a closed integral current
$C=C_2(S,\Theta)$ calibrated by $\Omega$, such that
cohomologically,
\begin{equation}
\lab{eq:chernofgasd} [C_2(A)] + {\rm PD}[C_2(S,\Theta)] = C_2(E),
\end{equation}
where $C_2(A)$ denotes the Chern-Weil form of $A$\footnote{One can
show this form, which was originally defined on $M\backslash
S(A)$, extends to a well-defined current on $M$.} and $C_2(E)$ denotes
the second Chern class of $E$. Two generalized $\Omega$-self-dual
instantons $(A, C)$, $(A', C')$ are equivalent if and only if
$C=C'$ and there is a gauge transformation $\sigma$ on
$M\backslash S(A) \cup S(A')$, such that $\sigma(A)= A'$ on
$M\backslash S(A) \cup S(A')$. We denote by $[A, C]$ the gauge
equivalence class of $(A,C)$. We identify $[A, 0]$ with $[A]$ in $
{\cal M}_{\Omega, E}$ if $A$ extends to a smooth connection of $E$
over $M$ modulo a gauge transformation. We define
$\overline{\M}_{\Omega, E}$ to be set of all gauge equivalence
classes of generalized $\Omega$-self-dual instantons of $E$ over
$M$.

The topology of $\overline{\M}_{\Omega, E}$ can be defined as
follows: a sequence $[A_i, C_i]$ converges to $[A,C]$ in
$\overline{\M}_{\Omega, E}$ if and only if there are
representatives $(A_i, C_i)$ such that their associated currents
$C_2(A_i, C_i)$ converge weakly to $C_2(A,C)$ as currents, where
\begin{equation}
\lab{eq:chernclassforgasd}
 C_2(A', C')= C_2(A') + C_2(S',\Theta'),
~~~C'=(S', \Theta').
\end{equation}
It is not hard to show that by taking
a subsequence if necessary, $\tau_i(A_i)$
converges to $A$ outside $S(A)$ and the support of $C$ for some
gauge transformations $\tau_i$.

The following was proved in \cite{tian5} and provides a compactification
for the moduli space of $\Omega$-self-dual connections.

\begin{theo}
\label{th:tianasd} For any $M$, $g$, $\Omega$ and $E$ as above,
${\overline \M}_{\Omega, E}$ is compact with respect to this
topology.
\end{theo}

Of course, ${\overline \M}_{\Omega, E}$ admits a natural
stratification. The remaining problem, which is also important for
issue $3$, is about regularity of a generalized $\Omega$-self-dual
instanton. Another interesting problem is to develop a deformation theory of
smoothing singular self-dual instantons. Are there any constraints
on a singular self-dual instanton which is the limit of smooth
self-dual instantons? We do not even know any example of
a Hermitian-Yang-Mills connection with an isolated singularity
and which can be approximated by smooth Hermitian-Yang-Mills connections.

Let us give an example. Assume that $(M,g)$ is a K\"ahler manifold
with K\"ahler form $\omega$. Put $\Om= \omega^{m-2}/ (m-2)!$,
where $n=2m$. Then an $\Om$-self-dual instanton $A$ is simply a
Hermitian-Yang-Mills connection, that is $F^{0,2}_A=0$ and
$F^{1,1}_A \cdot \omega =0$, where $F_A^{k,l}$ is the $(k,l)$-part
of $F_A$. If $(A,C)$ is a generalized $\Omega$-self-dual
instanton, it follows from a result of J. King that there are
positive integers $m_a$ and irreducible complex subvarieties $V_a$
such that for any smooth $\varphi$ with compact support in $M$,
$$C_2(S, \Theta) (\varphi) = \sum_a m_a \int_{V_a} \varphi.$$
On the other hand, using a result of Bando-Siu, one can show
\cite{tianyang} that there is a gauge transformation $\tau$ such
that $\tau(A)$ extends to be a smooth connection outside a complex
subvariety of codimension greater than $2$.

We expect that general self-dual connections have analogous
properties. If $(A,C)$ is a generalized $\Omega$-self-dual
connection, we would like to have (1) the regularity of the
current $C$, that is, $C$ is presented by finitely many
$\Omega$-calibrated subvarieties with integral multiplicity; (2)
There is a gauge transformation $\sigma$ such that $\sigma(A)$
extends to a smooth connection outside a subvariety of codimension
at least $6$. In \cite{uhlenbeck2}, any isolated singularity of a
Yang-Mills connection in dimension 4 can be removed.
In \cite{taotian}, a removable singularity theorem
was established for stationary Yang-Mills connections in higher dimensions.
Using this, we can conclude that $\sigma(A)$ extends to a smooth connection
outside a closed subset $S$ with vanishing $(n-4)$-Hausdorff
measure. Further understanding on $S$ is needed. We will discuss
regularity of $C$ in the next section.

A particularly interesting case is the complex self-dual instanton.
We do expect to construct new invariants for Calabi-Yau 4-folds by
showing that the above moduli space of generalized complex
self-dual instantons gives rise to a fundamental class. The main
problem left is the regularity of generalized instantons. A
special case of this can be done nicely. If the underlying
Calabi-Yau 4-fold $M$ is of the form $Y\times T_\CC^1$, where $Y$
is a Calabi-Yau 3-fold and $T_\CC^1$ is a complex 1-torus, then a
$T^1_\CC$-invariant complex self-dual instanton is given by a
Hermitian Yang-Mill connection $A$ on $Y$ and a $(0,3)$-form $f$
with $\overline{\partial}^*f=0$. The expected dimension of its moduli space
is zero. Counting them with sign gives rise to the holomorphic
Casson invariant, which was constructed previously by R. Thomas
using the virtual moduli cycle construction in algebraic geometry\cite{thomas}.

Other analytic problems on the Yang-Mills equation include whether
or not the Yang-Mills flow develops singularity at finite time. It
was proved by Donaldson that the Yang-Mills flow along Hermitian
metrics of a holomorphic bundle has global solution. Of course, if
the dimension of the underlying manifold is less than $4$, the
Yang-Mills flow has a global solution. In general, it is still
open. If a singularity forms at finite time, how does it look
like?

\Section{Minimal submanifolds}

\vskip-5mm \hspace{5mm }

The study of minimal submanifolds is a classical topic. We will
not intend to cover all aspects of this topic. We will only
discuss issues related to previous discussions and particularly
self-dual type solutions of the minimal submanifold equation.

Let $(M,g)$ be an $n$-dimensional Riemannian manifold and $S$ be a
submanifold in $M$. Recall that $S$ is minimal if its mean
curvature $H_S$ vanishes. The mean curvature arises from the first
variation of volume of submanifolds. Minimal submanifolds are
closely related to the Yang-Mills equation. In fact, it was shown
in \cite{tian5} that if a Yang-Mills connection has its curvature
concentrated along a submanifold, then this submanifold must be
minimal and of codimension 4. Motivated by this, recently, S.
Brendle, etc. developed a deformation theory of constructing
Yang-Mills connections from minimal submanifolds.

Now assume that $M$ has a closed differential form $\Omega$ with
its norm $|\Omega| \le 1$. A submanifold $S$ is calibrated by
$\Omega$ if $\Omega|_S $ coincides with the induced volume form.
Calibrated submanifolds are minimal (cf. \cite{harveylawson}). The
study of calibrated submanifolds was pioneered in the seminal work
\cite{harveylawson} of Harvey and Lawson. It now becomes extremely
important in the string theory. As we have seen in the above, they
also appear in formation of singularity in the Einstein equation.
In particular, when a self-dual connection has its curvature
concentrated along a submanifold, this submanifold is calibrated
\cite{tian5}, so a calibrated submanifold can be regarded as a
self-dual solution of the minimal submanifold equation.

Let $(M,\omega)$ be a symplectic manifold and $J$ be a compatible
almost complex structure, that is, $\omega(Ju, Jv)=\omega(u,v)$
and $\omega(u,Ju) > 0$ for any non-zero tangent vectors $u$ and
$v$. Define a compatible metric $g$ by $g(u,v)=\omega(u,Jv)$. Any
$\omega$-calibrated submanifolds are $J$-holomorphic curves, that
is, each tangent space is a $J$-invariant subspace in $TM$. They
are particularly minimal with respect to $g$. Holomorphic curves
have been used to establish a mathematical foundation of the
quantum cohomology, the mirror symmetry, etc. (cf.
\cite{ruantian}). The key of it is to construct the Gromov-Witten
invariants by showing the moduli space of $J$-holomorphic curves
can be taken as a fundamental class in a suitable sense. This was
proven by first constructing a ``nice'' compactification of the
moduli space of $J$-holomorphic curves and then applying
appropriate transversality theory.

I do believe that there should be new invariants by using other
calibrated submanifolds. A particularly interesting case is the
Cayley cycles in a ${\rm Spin}(7)$-manifold. Note that a
Calabi-Yau 4-fold is a special ${\rm Spin}(7)$-manifold. Again the
problem is about the structure of singular Cayley cycles. This
proposed new invariants will provide a powerful tool of
constructing Cayley cycles in a ${\rm Spin}(7)$-manifold,
particularly, holomorphic and special Lagrangian cycles in a
Calabi-Yau 4-fold. A related problem is to construct new
invariants for hyperK\"ahler manifolds by using tri-holomorphic
maps. A good compactification for the moduli of tri-holomorphic
maps is needed, but this should be technically easier. Partial
results have been obtained in \cite{litian} and \cite{chenli1}.

Another possible invariant may exist for Calabi-Yau manifolds. Let
$(M,\omega)$ be a Calabi-Yau $n$-fold with a holomorphic $n$-form
$\omega$ such that
\begin{equation}
\lab{eq:normalization}
\omega^n (-1)^{n(n-1)\over 2} n!
\left ({\sqrt{-1}\over 2}\right)^n \Omega\wedge{\overline {\Omega}}.
\end{equation}
A special Lagrangian submanifold is a submanifold $L\subset M$ such that
$\omega|_L=0$ and $\Omega$ restricts to the induced volume form of $L$.
If one has a good compactification theorem for special Lagrangian submanifolds,
then one can count them to obtain a new invariant for $M$. A particularly
important case is for Calabi-Yau 3-folds.

The minimal equation is nonlinear and does have singular
solutions. So one has to introduce weak solutions. An integral
$k$-dimensional current $C=(S, \Theta, \xi)$ consists of a
$k$-dimensional rectifiable set $S$ \footnote{This implies that
there is a unique tangent space at a.e. point of $S$.} of locally
finite Hausdorff measure, an $H^k$-integrable integer-valued function $\Theta$ and
a $k$-form $\xi \in \wedge^kTS$ with unit norm. Each current
induces a natural functional $\Phi_C$ on smooth forms with compact
support: For any smooth form $\varphi$,
\begin{equation}
\lab{eq:current} \Phi_C (\varphi) = \int_S \langle \varphi,
\xi\rangle dH,
\end{equation}
where $dH^k$ denotes the $k$-dimensional Hausdorff measure. We say $C$ has no
boundary if $\Phi_C(d\psi)=0$ for any $\psi$. One can define the
generalized mean curvature of $C$ as the variation of volume. A
current $C$ is minimal if its mean curvature vanishes. A current
$C$ is calibrated by a $k$-form $\Omega$ if $\Omega|_{T_xS}$
coincides with the induced volume form whenever the tangent space
$T_xS$ exists. Of course, a calibrated current is minimal,
provided that $d \Omega=0$ and $|\Omega|\le 1$.

A fundamental problem in the regularity theory of minimal surfaces
is the regularity of minimizing currents. A result of F. Almgren
claims that an area minimizing current is regular outside a subset
of Hausdorff codimension two \cite{almgren}. In many geometric
applications, we will encounter with calibrated currents, for
example, in the famous work of Taubes on equivalence of the
Seiberg-Witten invariants and the Gromov invariants, the key
technical point is to show that any $\omega$-calibrated current in
a symplectic 4-manifold $(M,\omega)$ is a classical minimal
surface, i.e., the image of a pseudo-holomorphic map from a smooth
Riemann surface (\cite{taubes}, also see
\cite{tristantian}).\footnote{ This is a special case of the main
result in \cite{schang}, which states that a $2$-dimensional area
minimizing is a classical minimal surface. But Chang's proof
relies on some hard techniques of \cite{almgren}, so a
self-contained proof is very desirable.} This current is obtained
as an adiabatic limit of curvature forms of solutions of deformed
Seiberg-Witten equations. The problems of this type should also
occur when we study the Calabi-Yau manifolds near large complex
limits. Of course, this regularity problem also appears in
compactifying moduli spaces of calibrated cycles.

Here is what we think should be true: If $C=(S,\Theta,\xi)$ is a
$k$-dimensional calibrated current, then $S$ can be stratified
into $\coprod_{i} S_i$ such that each stratum $S_i$ is a smooth
manifold of dimension $i$, which is at most $k-2$, and $\Theta$ is
constant on each stratum.

If the calibrating form $\Omega$ is $\omega^l/l!$ ($k=2l$) on a
symplectic manifold $(M,\omega)$ with a compatible metric $g$,
then $\Omega$-calibrated current is pseudo-holomorphic, that is,
any tangent space is invariant under the almost complex structure
induced by $\omega$ and $g$. In this special case, the conjecture
is that $S$ is stratified into pseudo-holomorphic strata $S_{2j}$.
If the dimension of $C$ is $2$, then the conjecture claims that
$C$ is induced by a pseudo-holomorphic curve. This conjecture
follows from a result of King when $(M,g)$ is
K\"ahler, i.e., the corresponding almost complex structure is
integrable. Very recently, Riviere and I can prove this conjecture
when $\dim C=2$. When $\dim M = 4$, it was already known (cf. \cite{taubes},
\cite{tristantian}).

A nice way of deforming a surface into a minimal one is to use the
mean curvature flow. If $(M,g)$ is a compact K\"ahler-Einstein
surface and $S_0$ is a symplectic surface with respect to the
K\"ahler form, then surfaces along the mean curvature flow
starting from $S_0$ are also symplectic \cite{chentian}. If the
flow has a global solution, then $S_0$ can be deformed to a
symplectic minimal surface. In particular, $S_0$ is isotopic to a
symplectic minimal surface. However, it is highly nontrivial to
show that the flow has a global solution. Partial results have
been obtained (\cite{chenli2}, \cite{mutaowang}). Nevertheless, it
was conjectured in \cite{tian6} that any symplectic surface in a
K\"ahler-Einstein surface is isotopic to a symplectic minimal
surface. This has been checked for a quite big class of symplectic
surfaces in a K\"ahler-Einstein surfaces with positive scalar
curvature (cf. \cite{sieberttian}) by using pseudo-holomorphic
curves. One can also ask similar questions for the mean curvature
flow along Lagrangian submanifolds. It was proved first by Smocsky
\cite{smoczyk}
that the mean curvature flow preserves the Lagrangian property. We
refer \cite{thomasyau} for more discussions on the mean curvature
flow for Lagrangian submanifolds.

\label{lastpage}

\end{document}